\input amstex
\documentstyle{amsppt}

\hsize=4.75in
\vsize=8in

\def\G{\Gamma}

\rightheadtext {Generation properties for von Neumann algebras}
\leftheadtext {Florin R\u adulescu}
\topmatter
\title  Finite generation type properties for fuchsian group
von Neumann algebras
\endtitle

\author Florin R\u adulescu
\footnote{
Research supported in part by the grant  DMS 9622911 from the National
Science Foundation} 
\footnote{Member of the Institute of Mathematics, 
Romanian Academy, Bucharest}\endauthor

\affil
Permanent Address: Department of Mathematics, The University of Iowa,
Iowa City, IA 52246, U.S.A 
\endaffil
\endtopmatter
\document
\bigskip

\def\G{\Gamma}
In this paper we prove the existence of  finite systems  $\Cal G$ of
generators, with special properties, for the von Neumann algebras
 $\Cal L(\G)\overline{ \otimes} B(H)$,
obtained by tensoring the group von Neumann algebra of
 a fuchsian group
$\G$ with the bounded linear operators on a infinite dimensional,
separable, Hilbert space $H$. The results show that
the linear span of ordered monomials that are products of 
powers of  adjoints of  elements in $\Cal G$ with
powers of  elements in $\Cal G$, is weakly dense
 in the algebra.
In addition the systems of generators are a  commuting family, having
non-trivial, joint invariant subspaces affiliated to the von Neumann
algebra.
 If $\G=PSL(2,\Bbb Z)$, we can find such a generating set consisting of a single element.

The  construction in this paper may  also be used to provide
examples for 
 Toeplitz operator with unbounded symbol with unexpected behaviour and 
 non-closability for multiplication operators, with bounded
symbol, in  Sobolev type Hilbert spaces associated to 
differential operators.

We recall some of the notations that were used in [Ra3].
We let $H_t$ be the Hilbert space of square summable,
analytic functions on $\Bbb H$, with respect to the 
measure $\text{d}\nu_t=(\text {Im\ }
z)^{t-2}\text {d}\overline {z} \text {d} z $. In [Pu], [Sa] it was 
proven that there exist a one parameter family of irreducible,
projective unitary representation of $PSL(2,\Bbb R)$ on $H_t$
that extends the analytic discrete series of representations for 
$PSL(2,\Bbb R)$. Moreover the methods in [VFR] along with the trace 
formula in [Pu], were used in ([Ra3]) to show that for any fuchsian 
group $\G$, the von Neumann algebra $\{\pi_t(\G\}''$ is a type 
$II$ factor acting on $H_t$. The Murray von Neumann dimension 
for this algebra acting on $H_t$ is proportional to $t$ and the
covolume of $\G$. In particular the commutant
algebras $\Cal A_t=\{\pi_t(\G\}"$ are isomorphic to
the twisted group von Neumann  algebra
 $\Cal L(\G,\sigma^t)_{(t-1)/\text{covolume \ }\G}$, with 
$\sigma^t$ the cocycle corresponding to the projective representation
$\pi_t$.

As shown in [GHJ], automorphic forms $f$ of weight $k$
for $\G$ correspond to bounded intertwining operators
$S^t_f$ from $H_t$ into $H_{t+k}$, by multiplication with $f$
on $H_t$. The intertwining property means that
$$\pi_{t+k}(\gamma) S^t_f=S^t_f\pi_t(\gamma), \gamma \in \G.$$
Let $f$ be any
measurable, bounded $\G$-invariant function on $\Bbb H$.
Let $M^t_f$ be the multiplication operator by $f$ on $H_t$
and let $P_t$ be the projection operator from bounded square summable
functions on $\Bbb H$ onto $H_t$. The 
Toeplitz operator $T^t_f=P_tM^t_fP_t$ clearly commutes with                     
$\{\pi_t(\G\}$ and hence $T^t_f$ belongs to $\Cal A_t$.

A natural object to consider is the (unbounded)
Toeplitz operator  with symbol a $\G$ invariant, analytic function on
$\Bbb H$. Unfortunately, as we  show bellow, though affiliated with
the commutant algebras, such an operator is not densely defined, nor is its
adjoint closable. One possible symbol function, for 
$\G=PSL(2,\Bbb Z)$, is the classical modular invariant function
  $j=\frac{\Delta}{G_4^3}$.

\proclaim {Lemma} Let $\G$ be a fuchsian group and let $j$ be 
an analytic, $\G$-invariant function on $\Bbb H$. Assume that
 $j=\frac{a}{b}$
is the quotient of two automorphic forms for $\G$ having the same
 weight $m$.
Assume that the functions $z\rightarrow \vert a(z)\vert^2
(\text{Im\ }z)^{m-2},$ $z\rightarrow \vert b(z)\vert^2
(\text{Im\ }z)^{m-2},$ are bounded. Also we assume that $a,b$ have not all
their zeroes and poles in common.

 For a bounded, measurable, $\G$ invariant
function $f$ on $\Bbb H$ we let $T^t_f$ be the Toeplitz operator
on $H_t$ with symbol $f$. Then, clearly, $T^t_f$ belongs to the
commutant
$\Cal A_t=\{\pi_t(PSL(2,\Bbb Z)\}'$ and the set
 of all such operators is a 
weakly dense  subspace of the Hilbert space $L^2(\Cal A_t)$ associated 
with the trace on $\Cal A_t$ (see [Ra1]).

   Then the operator 
$T^t_f\rightarrow T^t_{\overline {j} f}$ is not closable on
$L^2(\Cal A_t)$ for any $t>1$ (and hence
$T^t_f\rightarrow T^t_{ {j} f}$  is not densely defined). Moreover the
same holds true if we replace $\overline {j}$ by $\vert j\vert^2$ or
by $\frac{\overline {j}}{j}$.
\endproclaim 

Proof. The hypothesis shows that the bounded (see [GHJ])  operators $S^t_a, S^t_b$ 
on $H_t$ into $H_{t+k}$, defined by by multiplication with $a$ and 
respectively $b$ have non equal ranges. Let $c$ be any other 
automorphic form form $\G$ having the same order as $a,b$.
Let $\Psi$ be the operator in the statement. Then for any
$\G$-invariant, bounded, measurable function $h$ we will have 
that $$\Psi((S^t_a)^* T^{t+k}_hS^t_c)=(S^t_b)^* T^{t+k}_hS^t_c.$$
Let $e_a$ be the projections onto the range of $S_a^t$, which is a 
subspace of $H_{t+k}$.
As the set of all $T^{t+k}_h$, when $h$ runs trough the bounded,
 measurable,
 $\G$-equivariant 
functions on $\Bbb H$ is weakly dense in the commutant $\Cal A_{t+k}$
it follows that for any $x$ in $\Cal A_{t+k}$  we may
 find a sequence of such functions $h_n$ on $\Bbb H$
 such that $T^{t+k}_{h_n}$ converges weakly to 
$(1-e_a)x$ and hence $(S^t_a)^* T^{t+k}_{h_n}S^t_c$
converges weakly to 0.
 Then $(S^t_b)^* T^{t+k}_{h_n}S^t_c$ converges weakly
to $(S^t_b)^*(1-e_a)xS^t_c$. If $\Psi$ was closable  
  then it would follow that $(S^t_b)^*(1-e_a)xS^t_c$ is zero
 for all $x$ in $\Cal A_t$.
Since $x$ is arbitrary and since the   operators $S^t_a, S^t_b$
have non equal ranges, this is impossible. This completes the proof
for $\overline{j}$. The statement for the other two functions is
proved in a similar way:  multiplication by $\vert j\vert^2$ maps
$(S^t_a)^* T^{t+k}_hS^t_a$ into $(S^t_b)^* T^{t+k}_hS^t_b$ (and this
is a positive map) and multiplication by 
$\frac{\overline {j}}{j}$ maps $(S^t_a)^* T^{t+k}_hS^t_b$ into
$(S^t_b)^* T^{t+k}_hS^t_a$.
\bigskip

\proclaim{Remark}With the notations in the above statement
let $F$ be a fundamental domain for the action of $\G$ on $
\Bbb H$. Let $B_t=B_t(\Delta) $ be the Berezin operator
commuting with the invariant laplacian, formally defined
by 
 $$B_t f(z)=\int_{\Bbb H}f(w)\lbrack \frac{(\text {Im\ } z )
 (\text {Im\ } w)}{\overline {z}- w}\rbrack^{t-2}
(\text {Im\ } w)^{-2} \text {d} \overline{w}
 \text {d} w,\ \  z\in \Bbb H.$$

Then $L^2(\Cal A_t,\tau)$ is identified with the the completion
of a dense subset 
of\break $L^2(F,(\text {Im\ } w)^{-2} \text {d} \overline{w}
 \text {d} w) $ with the scalar product 
$$\langle f, g\rangle = (B_t(\Delta) f, g)_{
L^2(F)}.$$

With this scalar product, the operator of multiplication by 
$\overline {j}$
on a dense subspace of
 $L^2(F,(\text {Im\ } w)^{-2} \text {d} \overline{w}
 \text {d} w)$ is non closable. Moreover the
same holds true if we replace $\overline {j}$ by $\vert j\vert^2$ or
by $\frac{\overline {j}}{j}$.

\endproclaim

The following corollary shows that one of the properties that
 is 
valid for Toeplitz operators with bounded, antianalytic symbol, fails
for unbounded symbols (see also [Ja], [Saf]).

\proclaim{Corollary} Let $\G$ be a fuchsian group of finite covolume,
and $j=\frac{a}{b}$ an analytic, $\G$ invariant function on $\Bbb H$
as above. Let $P_t$ be the projection from $L^2(\Bbb H, (\text {Im\ }
z)^{t-2}\text {d}\overline {z} \text {d} z) $ onto 
$H_t=H^2(\Bbb H, (\text {Im\ }
z)^{t-2}\text {d}\overline {z} \text {d} z) $. Let $M_{\overline j}$
 be
the multiplication operator with $\overline {j}$ which is
defined  on a dense subset of the Hilbert space $L^2(\Bbb H, (\text {Im\ }
z)^{t-2}\text {d}\overline {z} \text {d} z) $.

Then $P_t M_{\overline j} (1-P_t)$ is nonzero (in particular has 
nonzero domain).

\endproclaim

Proof. Assume the contrary.
 Let $h,h_1$ be any  $\G$-invariant, bounded, measurable, real valued
 functions $h, h_1$ on $\Bbb H$, such that the operators
$T^t_h$, $T^t_{h_1}$ are injective.
Assume that the supports of $h,h_1$ are so that the functions
$\overline{j} h, \overline{j} h_1$ are bounded.

  Let 
$Z_h$, $Z_{h_1}$ be the closable operators ([MvN]), defined by 
$$Z_h=T^t_{\overline{j} h}(T^t_h)^{-1},
Z_{h_1}=T^t_{\overline{j} h_1}(T^t_{h_1})^{-1}.$$

Our assumption implies that
$$T^t_{\overline{j} h}\zeta =T^t_{\overline{j} h_1}\zeta_1
\ \text{whenever\ } T^t_h\zeta=T^t_{h_1}\zeta_1.$$

This implies that the closable, unbounded operators 
$Z_h$, $Z_{h_1}$ coincide on a densely defined core and hence that
they are equal.

Hence there exists a unique, closed operator, affiliated with
$\Cal A_t=\{\pi_t(\G)\}'$ such that for any real valued,
 bounded, measurable
functions $h$, with $T^t_h$ invertible, we have (by [MvN])
$$ZT^t_h=T^t_{\overline{j} h}.$$
But this is impossible by the previous statement.

\proclaim{Corollary} With the notations in the previous statement
the same conclusion holds if $j$ is replaced by $h\circ j$,
where $h$ is any univalent entire function $h$.
\endproclaim

Proof.
By examining the above argument, we see that
 in fact we  proved that 
if $K$ is any compact subset of the interior of $F$ such that
$j$ is bounded when restricted to $K$ and $\chi_K$ is the
characteristic function of $K$  then 
$$P_tM_{\overline j} \lbrack(\text{Id\ }-P_t) \land \chi_K\rbrack
\ne 0.$$

Now assume that the (bounded) linear operator
$M_{\overline{f\circ j}}$ would have the property that
$$M_{\overline{f\circ j}}[(\text{Id\ }-P_t) \land \chi_K)(H_t)]
\subseteq (\text{Id\ }-P_t)(H_t).$$
It would the follow that $M_{\overline{f\circ j}}$
 would also have the property that
$$M_{\overline{f\circ j}}[(\text{Id\ }-P_t) \land \chi_K)(H_t)]
\subseteq (\text{Id\ }-P_t) \land \chi_K)(H_t).$$
Since $f$ is univalent the same would then hold true about $j$
instead of $f\circ j$, and this we know to be false.

\proclaim{Questions}
\item (i). Is the above statement true if one drops the univalence
condition on $f$?
 
\item (ii). If $t>24$, $\G=PSL(2,\Bbb Z)$ then if 
$j=\frac{\Delta}{G^3_4}=\frac{\Delta^2}{G^3_4\Delta}$ it is clear
that that the domain of $M_{\overline j}$ intersects $H_t$ non trivially.
Does this hold for smaller values of $t$?
\endproclaim

Assume that $\G$ is a fuchsian group
(necessary of infinite covolume). Assume that
$H^{\infty}(\Bbb H/\G)$ has a sufficiently rich structure so that
there are a finite number of bounded analytic functions that
separate points on $\Bbb H/\G$. In this case, by using the methods
developed in [Ra3] we can show that $\Cal L(\G)
\otimes B(K)$ has a  set of generators with  the properties in the following
proposition.

\proclaim{Proposition 2} Let $\G$ be a fuchsian group such that
 $H^{\infty}(\Bbb H/\G)$ contains functions $h_1, ... h_k$ that
separate the points on $\Bbb H$. Let $K$ be an infinite dimensional 
Hilbert space. Then there exists commuting,
 bounded operators $Z_1,...Z_k$
in $\Cal B= \Cal L(\G)\otimes B(K)$ such that 
$$\Cal B= \overline {
\text{Sp\ } \{ (Z_i^n)^\ast Z_j^m\vert
n,m \in \Bbb N, i,j=1,2...,k\}}^{\text{weak}}.$$
\endproclaim

Proof. Let $Z_i=T^t_{h_i}, i=1,...,k$. Let $c$ be an arbitrary
invertible element in $\Cal B\cap L^1(\Cal B,\tau)$. Assume that
$a\in \Cal B$ is orthogonal to the following subspace of
 $L^1(\Cal B,\tau)$
$$\overline{\text{Sp\ }\{ c(Z^n_i)^\ast Z^m_j \vert
n,m\in \Bbb N,i=1,...k\}}^{\text {weak}}.$$
We use the notations in [Ra3].
 Let $\hat{(ac)}(\overline{z},z), z\in \Bbb H$
be the Berezin's symbol of $ac \in \Cal B\cap L^1(\Cal B,\tau)$.
 The trace formula in [Ra3] shows that 
$$\int _F \hat{(ac)}(\overline{z},z)
\overline {h_i^n(z)}h_j(z)(\text {Im\ }
z)^{t-2}\text {d}\overline {z} \text {d} z=0.$$
Hence $ac=0$ and hence $a$ =0. Consequently
$$\text{Sp\ }\{ c(Z^n_i)^\ast Z^m_j \vert
n,m\in \Bbb N,i=1,...k\}$$
is a weakly dense subspace of $L^1(\Cal B,\tau)$ and consequently,
since $c$ was invertible we get that
$$\overline{\text{Sp\ }\{ (Z^n_i)^\ast Z^m_j \vert
n,m\in \Bbb N,\ i,j=1,...,k\}}^{\text {weak}}=\Cal B.$$

\proclaim{Proposition} Let $E$ be an open $\G$-invariant subset of $\Bbb H$, such that $j$
(or $\frac{1}{j}$) is bounded on $E$. Let $H_t(E)$ be the
 subspace of $L^2(\Bbb H,\nu_t)$, ($\G$-invariant, with respect to the unitary representation
of $PSL(2, \Bbb R)$ on $L^2(\Bbb H,\nu_t)$), consisting of square integrable, analytic
functions on $E$ (that are extended with 0 outside $E$). Let $Z$ be the Toeplitz operator on $H_t(E)$ with symbol
$j|_E$. Then $H_t(E)$ is a (finite or infinite) Hilbert module over $\Cal L(\G)$, as a submodule of $L^2(\Bbb H,\nu_t)$.

Then the commutant of $\G$ in $B(H_T(E))$ is the  weak closure of the 
linear span of the set $\{(Z^*)^n Z^m| n,m=0,1,2...\}$.
\endproclaim

Proof. Clearly $H_t(E)$ is a submodule of  $L^2(\Bbb H,\nu_t)$ over $\Cal L(\G)$. If $k_E$ is the reproducing kernel
for $H_t(E)$ and $F(E)$ is a fundamental domain for $\G$ acting on $E$, then the Muray von Neumann dimension of
$H_t(E)$ over  $\Cal L(\G)$ is equal to $\int_{F(E)}k_E(z,z)$. Moreover, as in [Ra3], if $A$ is an operator acting on $H_t(E)$,
we consider its Berezin kernel (with respect $H_t(E)$) to be 
$$k^E_A(z,\zeta)=\frac{<Ae^{E,t}_z,e^{E,t}_\zeta>}{<e^{E,t}_z,e^{E,t}_\zeta>},z,\zeta\in E,$$
where $e^{E,t}_z$, for $z$ in $E$ is the evaluation vector in$H_t(E)$ at $z$. Then , if $A$ is trace class in
the commutant of $\G$, the trace is $\int_{F(E)}k^E_A(z,z)\text{d}\nu_0(z)$.
 This trace comes from the trace on the commutant of $\G$ on $L^2(\Bbb H,\nu_t)$ that is normalized by giving value 
$\frac{t-1}{\pi}$ to $H^2(\Bbb H, \nu_t)$. As in [Ra3], if an operator is orthogonal on all Toeplitz operators
on $H_t(E)$, with $\G$-equivariant symbol it follows that $A=0$. Hence the linear span of Toeplitz operators,
with $\G$-equivariant symbols, is weakly dense in the commutant and since the symbol of $(Z^*)^n Z^m$ is 
$(\overline j|_E)^n(j|_E)^m$ the statement follows.

\proclaim{Corollary} In the algebra $\Cal A=\Cal L(PSL(2,\Bbb Z)\otimes B(H)=\Cal L(F_N)\otimes B(H)$, $N$ finite,
 there exists a bounded subnormal operator $Z$, such that $\Cal A$ is the weak closure of linear span of the set $\{(Z^*)^n Z^m| n,m=0,1,2...\}$.
\endproclaim

\centerline{References}

\item {[Be]} F. A. Berezin, Quantization in complex
  symmetric spaces,
Math USSR\break Izvestija, 9(1975),  341-379.

\item{[Co]} J. B. Conway, {\it Subnormal Operators}, Research Notes in Mathematics v. 51, 
Pitman Advanced Publ. Prog., Boston, 1982.

\item{[Dyk]}       K. Dykema, Free products of
 hyperfinite von Neumann algebras and
free dimension, Duke Math. J. 69/1 97-119 (1993).

\item{[GHJ]}       F. Goodman, P. de la Harpe, V.F.R. Jones, Coxeter
Graphs and Towers of Algebras, Springer Verlag, New York, Berlin,
Heidelberg, 1989.

\item{[JS]}  Janas, Jan
                 Stochel, J ,       
  Unbounded Toeplitz operators in the Segal-Bargmann space. II.

\item{[MvN]}, F. J. Murray, J. von Neumann, On ring of Operators,IV,
 Annals of Mathematics, 44 (1943), 716-808.

\item{[Pu]}      L. Pukanszki, The Plancherel formula
 for the universal covering group of
PSL(2, R), Math Annalen, 156 (1964), 96-143.

\item{[Pi]}     G. Pisier, Espaces de Banach quantiques:
 une introduction \` a la 
th\' eorie des espaces des operateurs, JournŽe Annuelle,
 Soc. Math. France, 1994.

\item{[Ra1]}       F, R\u adulescu,On the von 
Neumann Algebra of Toeplitz Operators with
Automorphic Symbol, in Subfactors, 
Proceedings of the Taniguchi Symposium on Operator
Algebras, edts. H. Araki, Y. Kawahigashi, H. Kosaki,  World Scientific,
 Singapore-New Jersey,1994.

\item{[Ra2]}      F. R\u adulescu, Random matrices,
 amalgamated free products and
subfactors in free group factors of noninteger index,
 Inv. Math. 115,  pp. 347-389, (1994).

\item{[Ra3].}F. R\u adulescu, The arithmetic Hecke operators
and Berezin Quantization,
{\it  Comptes Rendu Acad. Sci. Paris, 1996, Serie I.
                   MathŽmatique}
                 {\bf 322} (1996), no. 6, 541--546.

\item{[Ra4]}F. R\u adulescu, The $\Gamma$ invariant form of
 the Berezin quantization of the upper
halfplane (Preprint 1995, to appear in Memoirs A. M. S. ).

\item{[Sa]}     P. Sally, Analytic
 Continuation of the Irreducible
Unitary Representations of the Universal
 Covering Group, Memoirs A. M. S. , 1968.

\item{[St]} Stout, E. L.,
 On some algebras of analytic functions on finite open Riemann
                   surfaces.
 Math. Z., 92 1966 366--379.

\item{[Vo]}       D. Voiculescu, Circular and semicircular
 systems and free product
factors. In Operator Algebras, Unitary Representations,
 Enveloping algebras and Invariant Theory
(Prog. Math., v. 92, pp. 45-60, Boston: Birkhauser 1990.

     \end